\documentclass{notices}
\usepackage{amsfonts, amssymb, amsmath, amscd}
\usepackage{esint}
\usepackage{amsthm}
\usepackage[dvipsnames]{xcolor} 
\usepackage{enumerate}
\usepackage{mathtools}
\usepackage{enumitem}
\usepackage{bbm}
\usepackage{microtype}
\usepackage{hyperref}
\usepackage{tikz}
\usepackage{pgfplots}
\pgfplotsset{compat=1.15}
\usepgfplotslibrary{fillbetween}
\usetikzlibrary{patterns, angles, quotes}
\usepackage{tikz-3dplot}
\usepackage{url}
\usepackage{dsfont}
\usepackage{longtable}

\def\eqn#1$$#2$${\begin{equation}\label#1#2\end{equation}}

\def\loc{\operatorname{loc}}

\numberwithin{equation}{section}
\allowdisplaybreaks[4]

\newcommand{\tta}{s}

\newcommand{\hhh}{\textnormal{\texttt{h}}}
\newcommand{\mmm}{m}

\newcommand{\kk}{\kappa}

\def\er{\mathbb R}

\newcommand{\ti}[1]{\tilde{#1}}
\newcommand{\mf}[1]{\mathfrak{#1}}

\newcommand\eps\varepsilon

\def\eqn#1$$#2$${\begin{equation}\label#1#2\end{equation}}

\newcommand{\be}{\begin{equation}}
\newcommand{\ee}{\end{equation}}

\def\BBB{\textnormal{B}}

\newcommand{\rr}{\varrho}

\newcommand{\snr}[1]{\lvert #1\rvert}

\newcommand{\nr}[1]{\lVert #1 \rVert}

\newcommand{\avnorm}[1]{\linethrough\|#1 \|}

\def\name[#1, #2]{#1 #2} 

\newcommand{\sss}{\mf{s}}

\newcommand{\rif}[1]{(\ref{#1})}

\newcommand{\medint}{-\kern -,375cm\int}

\definecolor{gialloparma}{rgb}{0.99,0.845,0}

\definecolor{bluparma}{rgb}{0.016,0.51,0.79}
\newcommand{\medintinrigo}{-\kern -,315cm\int}
\makeatletter
\newcommand{\linethrough}{\mathpalette\@thickbar}
\newcommand{\@thickbar}[2]{{#1\mkern0mu\vbox{
    \sbox\z@{$#1#2\mkern-0.5mu$}%
    \dimen@=\dimexpr\ht\tw@-\ht\z@+2\p@\relax %
    \hrule\@height0.5\p@ 
    \vskip\dimen@
    \box\z@}}
}
\makeatother

\def\aaa{\mf{a}}

\newcommand\ccc{\mathfrak{c}}

\def\dx{\,{\rm d}x}

\def\dy{\,{\rm d}y}

\def \d{\,{\rm d}}
\def \diver{\,{\rm div}}
\def\dist{\,{\rm dist}}

\def\kkk{{\textnormal{\texttt{k}}}}

\def\deb{\rightharpoonup}

\theoremstyle{plain}
\newtheorem{thm}{Theorem}

\newtheorem{cor}[thm]{Corollary}

\newcommand{\thistheoremnames}{}
\newtheorem*{genericthms}{\thistheoremnames}
\newenvironment{para*}[1]
  {\renewcommand{\thistheoremnames}{#1}%
   \begin{genericthms}}
  {\end{genericthms}}

\theoremstyle{remark}

\numberwithin{equation}{section}

		\def \diver{\,{\rm div}}

\title{Nonuniform ellipticity in variational problems and regularity}

\author{Cristiana De Filippis and Giuseppe Mingione}

\allowdisplaybreaks

 \pgfplotsset{unit circle/.style={width=8cm,height=8cm,axis line style={draw=none},xtick=\empty,ytick=\empty,axis equal,enlargelimits,xmax=1,ymax=1,xmin=-1,ymin=-1,domain=-1:1}}

\begin{document}

\maketitle

\setcounter{tocdepth}{1}

\section{Existence, regularity}\label{zerosec}
Since Hilbert's 19 and 20th problems, the regularity theory of elliptic equations has undergone extensive development over the past century, often in connection with the minimization of integral functionals of the Calculus of Variations of the type  
\eqn{basicF}
$$
 \mathcal{F}(w,\Omega):= \int_{\Omega} F(x, Dw)\dx \,,\quad \Omega \subset \er^n
$$
defined for functions $w\in W^{1,1}(\Omega)$. 
Here $\Omega$ denotes a bounded open set, $n\geq 2$, and $F\colon \Omega \times \er^n\to [0, \infty)$ is such that $x \mapsto F(x,z)$ is measurable for every $z \in \er^n$ and $z \mapsto F(x,z)$ is convex for every $x \in \Omega$\footnote{As a consequence,  $x\mapsto F(x,Dw(x))$ is measurable.}. This in fact makes $\mathcal F$ itself a convex functional. In \eqref{basicF}  $Dw$ denotes the (distributional) gradient of $w$, while the Sobolev space $W^{1,p}(\Omega)$, $1 \leq p \leq \infty$, consists of all $L^p(\Omega)$-functions $w$ such that $Dw\in L^{p}(\Omega;\er^n)$. This is a Banach space when equipped with norm $\nr{w}_{W^{1,p}(\Omega)}:=\nr{w}_{L^p(\Omega)}+\nr{Dw}_{L^p(\Omega)}$
. Sobolev spaces provide a natural framework for the modern Calculus of Variations, striking an effective balance between regularity features of functions and their compactness properties in certain weak topologies. This balance allows to use the so-called Direct Methods in the Calculus of Variations - an infinite-dimensional analogue of Weierstrass's theorem. The classical, one-variable Weierstrass's theorem, states that a function $f\colon X (\subset \er) \to \er$ attains its minimum on $X$ provided $f$ is lower semicontinuous\footnote{that is, $ f(x) \leq \liminf_{i}  f(x_i)$ holds whenever $x_i\to x$.}, and $X$ is compact. It now occurs that, for $1 < p < \infty$, the space $W^{1,p}(\Omega)$ can be equipped with a (weak) topology making its bounded subsets sequentially precompact. The twist is that the convexity of $z\mapsto F(\cdot, z)$ implies the weak  sequential lower semicontinuity of $\mathcal F$, i.e., 
$
 \mathcal{F}(w,\Omega) \leq \liminf_{i}  \mathcal{F}(w_i,\Omega)
$
holds whenever a sequence $\{w_i\}\subset W^{1,p}(\Omega)$ converges to a function $w$ in the weak topology of $W^{1,p}(\Omega)$. 
At this point, we have all the necessary ingredients to revisit the original proof of Weierstrass’s theorem. Various existence results for minimizers can now be established - for example, by fixing boundary data - provided we can show the existence of at least one bounded minimizing sequence\footnote{In $W^{1,p}(\Omega)$ this can be reached for instance imposing a coercivity condition like $\snr{z}^p \lesssim F(x, z)$. A minimizing sequence is a sequence along which a functional attains its $\inf$.}. In short, these existence theorems give rise to the broad and unifying concept of local minimizer, which captures the essential minimality properties needed for the (interior) qualitative analysis of solutions and to which we reduce the analysis from now on. A function $u \in W^{1,1}_{\loc}(\Omega)$ is a local minimizer - or, from now on, simply a minimizer - if for every ball $\BBB\Subset \Omega$,  $F(\cdot, Du) \in L^1({\BBB})$ and $\mathcal{F}(u,\BBB)\leq \mathcal{F}(w,{\BBB})$ holds for every competitor $w \in u + W^{1,1}_0({\BBB})$\footnote{$W^{1,p}_0(\Omega)$, $1\leq p < \infty$, is the closure of smooth functions with compact support in $\Omega$ with respect to the $W^{1,p}$-norm. Elements of $W^{1,p}_0(\Omega)$ can be thought as Sobolev functions having ``zero boundary values".}. But the one-dimensional parallel can be pushed  one step further. Consider a minimizer $u$ of $\mathcal F$, a ball $\BBB\Subset \Omega$ and $\varphi\in C^{\infty}_0(\BBB)$, and observe that the one variable function $g(t):= \mathcal F(u+t\varphi, \BBB)$ reaches its minimum at $0$. Therefore, if $g'(0)$ exists, then $g'(0)=0$. In fact, under suitable technical assumptions\footnote{For instance, the following are sufficient: $u \in W^{1,p}(\BBB)$, $F(x,z)\lesssim \snr{z}^p+1$ and $z \mapsto F(\cdot, z)$ is $C^1$-regular.}, $g'(0)$ exists and can be computed as 
\eqn{debole}
$$
 \int_{\BBB} \partial_{z_{i}}F(x, Du)D_{i}\varphi\dx=0
$$
for every test function $\varphi\in C^{\infty}_0(\BBB)$  and every ball $\BBB\subset \Omega$ (Einstein notation over repeated indices is used here and throughout the paper). 
This is called the Euler-Lagrange equation of the functional $\mathcal F$ and, formally integrating by parts, can be written in {\em divergence form} as
\eqn{el0}
$$
-\diver\, \partial_zF(x, Du)=0\,.
$$ 
The catch with ellipticity is now that convexity of $z \mapsto F(\cdot, z)$ essentially imply that \rif{el0} is an elliptic equation - see next section. 

Regularity theory aims to establish smoothness properties of minimizers beyond those implied by their initial Sobolev function class, such as continuity, differentiability and higher-order smoothness. This is a fundamental issue for several reasons. For instance: minimizers and solutions to elliptic equations often represent certain relevant physical quantities in mathematical modelling. Since explicit solutions are generally unavailable, numerical schemes are developed to approximate them. Knowing a priori the regularity of solutions is essential for designing efficient approximation procedures.  Moreover, regularity is closely related to theoretical aspects. Establishing their regularity - e.g., proving $C^1$-regularity - demonstrates the existence of classical (strong) solutions, bridging the gap between the modern weak formulations seen above and more classical concepts of solution.

Before proceeding, some minimal notation. With $x\in \er^n$, $r>0$, we denote $B_{r}(x):=\{y \in \er^n \, \colon \, \snr{y-x}<r\}$. When the context is clear, we may simply write $\BBB\equiv B_{r}\equiv B_{r}(x)$; moreover, for $d>0$, we shall also denote $\BBB/d\equiv B_{r/d}(x)$ and use $\BBB\subset \er^n$ to denote a generic ball, regardless its center or radius. We shall also use averaged norms on balls $\BBB$, defined as 
$
\avnorm{f}_{L^{\kk}(\BBB)}:=\snr{\BBB}^{-1/\kk}\nr{f}_{L^{\kk}(\BBB)}
$
whenever $\kk>0$. Throughout the paper $c, \tilde c, c_1$ and similar symbols will denote generic constants $\geq 1$; writing $c\equiv c (t,s)$ indicates that $c$ 
depends on the parameters $t,s$. As usual, $u \lesssim_{t,s} v$ means $ u \leq cv$ for some $c\equiv c(s,t) $. Finally, for any $v\in \er$, we denote $v_+:= \max\{v,0\}$. 
 \section{Uniform, nonuniform}\label{prima}
Assuming that $z\mapsto F(\cdot, z)$ is $C^2$-regular outside the origin, when $\snr{z}\not=0$ we denote by $\lambda(x,z)$ and $\Lambda(x,z)$ the smallest and largest eigenvalues of $\partial_{zz}F(x,z)$, respectively, and we further assume $\lambda(x,z)>0$, that is, ellipticity. A crucial number quantifying the ellipticity properties of the integrand $F$ is its so-called ellipticity ratio $\mathcal R_{F}$, which is defined by 
$$
\mathcal R_{F}(x,z) := \frac{\Lambda(x,z)}{\lambda(x,z)}
$$  
for every $(x,z)\in \Omega \times(\er^n\setminus\{0_{\er^n}\})$. The integrand $F$ (or, equivalently, the functional  $\mathcal{F}$ and the equation \rif{el0}) is said to be uniformly elliptic if $\mathcal R_{F}$ is bounded in $\Omega \times (\er^n\setminus\{0_{\er^n}\})$. Otherwise, it is referred to as nonuniformly elliptic. 
In the following, we are interested in cases where the ellipticity ratio blows up as $\snr{z}\to \infty$.\footnote{Some authors refer to this as the notion of nonuniform ellipticity; in other words, the size of the ellipticity ratio is tested only for large values of  $\snr{z}$; see for instance \cite{simon}.} Classical examples of non-uniformly elliptic equations involve the minimal surface operator $\diver \, ((1+\snr{Du}^2)^{-1/2}Du)$ stemming from the choice $F(Dw)=(1+\snr{Dw}^2)^{1/2}$; in that case it is $\mathcal R_{F}(z)= \snr{z}^2+1$. Further examples include 
$$
\scalebox{0.91}{ $\displaystyle F(x,Dw)\equiv (|Dw|^2+1)^{p/2} + \sum_{k=1}^n \aaa_k(x)(|D_{k}w|^2+1)^{p_{k}/2}$}
$$ with $1 < p < p_{k} $, $0 \leq \aaa_k \leq L$ and   
$$
\begin{cases}
F(x,Dw) \equiv  \exp(\mf{c}(x)|Dw|^p),\quad  p\geq 1\\
 F(x,Dw) \equiv  \mf{c}(x)\snr{Dw}\log(\snr{Dw}+1)
 \end{cases}
 $$ with $1\leq \mf{c}\in L^{\infty}(\Omega)$. Here we focus on integrands with polynomial growth and H\"older continuous dependence on coefficients. 
Specifically, we adopt the following assumptions:
\eqn{assf}
$$
\scalebox{0.89}{$
\begin{cases}
\, \tilde{z}\mapsto F(\cdot,\tilde{z})\in C^{2}_{\loc}(\mathbb{R}^{n}\setminus \{0_{\er^n}\})\cap C^{1}_{\loc}(\mathbb{R}^{n})\vspace{2mm}\\
\, H_{\mu}(z)^{p/2}\le F(x,z)\le LH_{\mu}(z)^{q/2}+LH_{\mu}(z)^{p/2}\vspace{2mm}\\
\, H_{\mu}(z)^{(p-2)/2}\snr{\xi}\le  \partial_{z_iz_j}F(x,z)\xi_i\xi_j \vspace{2mm}\\
\, \snr{\partial_{zz}F(x,z)}\le L H_{\mu}(z)^{(q-2)/2}+LH_{\mu}(z)^{(p-2)/2}\vspace{2mm}\\
\, \snr{\partial_{z} F(x,z)-\partial_{z}F(y,z)}\vspace{2mm}\\ \, \quad \le L \snr{x -y}^{\alpha}[H_{\mu}(z)^{(q-1)/2}+H_{\mu}(z)^{(p-1)/2}]
\end{cases}$}
$$
for all $z, \xi \in \mathbb{R}^{n}$, $|z|\not=0$, $x,y\in \Omega$, where $p,q,\mu, L, \alpha$ are numbers such that $1<p \leq q $, $0 \leq \mu\leq 1 \leq L$, $0<\alpha \leq 1$, and where $H_{\mu}(z):=\snr{z}^{2}+\mu^{2}$ for any $z\in \er^n$. Assumption \rif{assf}$_3$ guarantees the strict convexity of $\mathcal{F}$, while \rif{assf}$_2$ ensures that any local minimizer of $\mathcal{F}$ automatically belongs to $W^{1,p}_{\loc}(\Omega)$. The parameter $\mu$ distinguishes the degenerate case $\mu=0$ from the non-degenerate one $\mu>0$. For example, the degenerate but uniformly elliptic integrand 
 $F(z)=\snr{z}^p$, $p>1$, satisfies \rif{assf} with $p=q$ and $\mu=0$. Its Euler-Lagrange equation is the familiar $p$-Laplace  equation 
$\diver (\snr{Du}^{p-2}Du)=0.$ The conditions in \eqref{assf} are called 
$(p,q)$-growth conditions, or non-standard growth conditions, highlighting the departure from the classical case $p=q$. The convexity and growth conditions \rif{assf}$_{3,4}$ allow to bound the ellipticity ratio as follows
\eqn{rate}
$$ 
\mathcal R_{F}(x,z) \lesssim_{p,q} L (\snr{z}^{q-p} +1)
$$ 
so that uniform ellipticity is ensured only when $p=q$.
In this situation it is clear that local $L^\infty$-gradient estimates are the focal point of regularity. Indeed, in Leon Simon's words \cite{simon}: ``Generally speaking, one can think of a local gradient bound as locally reducing a non-uniformly elliptic equation to a uniformly elliptic one, thus making it possible to use the theory of uniformly elliptic equations to study the solutions to non-uniformly elliptic equations."

\section{Gradient bounds}\label{gradb}
The regularity theory for nonuniformly elliptic problems is vast and, to a large extent, classical. This field has been systematically developed over several decades by numerous researchers, including  Bombieri, De Giorgi, Giusti, Gilbarg, Ivanov, Ivochkina, Ladyzhenskaya, Lieberman, Miranda, Serrin, Leon Simon, Stampacchia, Trudinger, Uraltseva, just to name a few. Providing a reasonable overview here is truly impossible. Therefore, we will limit our focus to the variational framework of Section \ref{prima}. The central strategy for establishing the regularity of minimizers involves  controlling the growth of the ellipticity ratio $\mathcal R_{\mathcal{F}}$ with respect to the gradient. In fact, by looking at the uniformly elliptic case, one sees that the boundedness of $\mathcal R_{\mathcal{F}}$ is always implicitly used when deriving gradient bounds. In nonuniformly elliptic problems a competition emerges: on one hand, one wants to exclude the possibility that the gradient of solutions blows up, while, on the other, gradient bounds rely critically on constraints imposed on the growth of $\mathcal R_{\mathcal{F}}$ with respect to the gradient. In this respect, the condition playing a central role is of the form 
\eqn{central}
$$
\frac qp < 1 +{\rm o}_n,
$$  
where ${\rm o}_n$ is a positive quantity such that ${\rm o}_n\to 0$ when $n\to \infty$. This implies that $q$ and $p$ cannot differ too much, ensuring that $\mathcal R_{\mathcal{F}}$ does not grow too rapidly, as inferred from \eqref{rate}. Giaquinta and Marcellini showed that if \rif{central} is not met for sufficiently small ${\rm o}_n$, then unbounded minimizers exist (see for instance \cite{paolo2}). On the positive side we have:
\begin{thm}\label{mainaut}
Let $u\in W^{1,1}_{\loc} (\Omega)$ be a minimizer of the functional $ \mathcal{F}$ in \eqref{basicF} with $F(x, z)\equiv F(z)$, under assumptions \eqref{assf}$_{1-4}$ and 
$$
\frac qp  < 1 + \frac2{n-1}.
$$  Then $Du$ is locally bounded in $\Omega$ and
 \eqn{milano2}
$$
\|Du\|_{L^{\infty}(\BBB/2)} \lesssim_{n,p,q,L}  \,   \avnorm{F(Du)}_{L^{1}(\BBB)}^{\mf{s}/q} + 1 
$$
holds for every ball $\BBB\Subset \Omega$, where  
 \eqn{nuovino}
$$ 
\scalebox{0.88}{$
\mf{s}:=\begin{cases}
\, \frac{2q}{(n+1)p-(n-1)q}\qquad \mbox{if either $n\ge 4 $ or $n=2$}\\[5pt]
\, \mbox{any number} > \frac{q}{2p-q}\quad \mbox{if $n=3$ and $q>p$}\\[5pt]
\, 1 \qquad \mbox{if $n=3$ and $q=p$}.
\end{cases}$}
$$
\end{thm}
Note that $\mf{s}\geq 1$ and $\mf{s}=1$ only when $p=q$. In this last case \rif{milano2} recovers the usual $L^p-L^\infty$ gradient estimate of the uniformly elliptic case (see \rif{reversa2} below). Moreover, $\mf{s}\to \infty$ when $q/p\to 1+2/(n-1)$. Theorem \ref{mainaut} originally appeared in Marcellini's pioneering paper \cite{paolo2} where \rif{milano2} was  proved under the more restrictive bound 
$q/p < 1 + 2/n$ and a slightly larger $\sss$. The version here, due to Bella and Sch\"affner (see for instance \cite{bella}), incorporates the sharpest bound on $q/p$ currently available in the literature under assumptions \rif{assf}. The optimal bound remains unknown, but progress has been made in \cite{cris2}, which features sharp results for $p=2$. The emphasis we place on the identity of $\sss$ in \rif{nuovino} is intentional; its precise value will play a role in Section \ref{stimesec}. 

A natural question arises regarding gradient estimates for nonautonomous integrals of the form \eqref{basicF}.  The methods available for Theorem \ref{mainaut} cannot be extended to this case unless $x\mapsto \partial_z F(x, \cdot)$ is differentiable and satisfies additional assumptions.  This stems from the fact that the techniques in \cite{paolo2, bella} and related works rely on differentiating the Euler-Lagrange equation \eqref{el0}, an approach precluded by \eqref{assf}$_5$. See Section \ref{inside} for further details. For uniformly elliptic equations, the treatment of H\"older coefficients relies on classical Schauder estimates, pioneered in the 1930s by Hopf, Giraud, Schauder, Caccioppoli. Referring to the linear model case:
\eqn{eq1}
$$
\begin{cases}
\, -\diver (a(x)Dv)=0  \quad \mbox{in $\Omega$}
\\
\,   \snr{\xi}^2 \leq  a_{ij}(x)\xi_i\xi_j\,, \quad  \nr{a}_{L^{\infty}(\Omega)}\leq L
\end{cases}
$$
for a.e. $x\in \Omega$ and every  $\xi \in \er^n$, one has the optimal  
\eqn{sch}
$$
a_{ij}  \in C^{0, \alpha} \Longrightarrow Dv \in C^{0, \alpha} \,, \quad 0<\alpha <1
$$ 
which holds for $W^{1,2}$-solutions $v$. 
Both local and global versions of \rif{sch} are available.
The proofs are based on perturbation arguments, that is, on a basic principle ruling elliptic equations: stability of solutions with respect to variations of data (in this case, coefficients). This means comparing $v$ on small balls $B_{\rr}(x_{0})$ to solutions $w$ of equations with ``frozen", i.e., constant coefficients $ \diver (a(x_{0})Dw)=0$ such that $w=v$ on $\partial B_{\rr}(x_{0})$. These solutions enjoy good regularity estimates, which are inherited by $v$, eventually leading to \rif{sch}. Several more recent and direct proofs of Schauder estimates for \rif{eq1} are known, avoiding the use of Potential Theory employed in the original approaches. For instance, Campanato made use of suitable function spaces characterizing H\"older continuity via mean square deviations of functions (integral oscillations). Trudinger gave a proof using certain decay properties of convolutions. Leon Simon's proof instead uses blow-up methods. In every case, the perturbative nature of the arguments employed remains unchanged. The perturbation approach extends to nonlinear equations, as shown by Giaquinta and Giusti \cite{giagiu}, who proved Schauder estimates for uniformly elliptic equations of the type $\diver A(x,Du)=0$ as \rif{el0} in the non-degenerate case $p=2$. The degenerate case is instead a classical result of Manfredi \cite{juan}, who for instance proved that $W^{1,p}$-solutions to $\diver ( \ccc(x) \snr{Du}^{p-2}Du)=0$, $1 \leq \ccc \leq L$, $p>1$, have locally H\"older continuous gradient provided $\ccc$ is locally H\"older continuous. Today, the literature on Schauder estimates is gigantic, encompassing results for a wide range of operators and settings, all fundamentally relying on perturbation methods. These estimates play a crucial role in establishing higher-order regularity of solutions via bootstrap techniques. A key point of interest here is the reliance of these perturbation methods on homogeneous a priori estimates. Such estimates naturally fail in anisotropic problems under conditions like \rif{assf}, unless  
$p=q$. This limitation has left the validity of Schauder estimates in nonuniformly elliptic settings an open question for a long time. This issue is not merely technical. As shown in the next section, Schauder estimates do not always hold in nonuniformly elliptic problems. Counterintuitively, in such a setting, coefficients are not ``perturbative".

\section{Disturbing examples}\label{lavsec} 
The functional
\eqn{doppio}
$$
\mathcal D_{\aaa}(w, \Omega) := \int_{\Omega} (\snr{Dw}^p+\mf{a}(x)\snr{Dw}^q)\dx\,,
$$
where $1<p <q$ and $0\leq \mf{a}\in L^{\infty}(\Omega)$, was first considered in the pioneering work of Zhikov \cite{vassili}. This is now widely known as the Double Phase functional and has inspired an extensive body of literature. The term ``Double Phase" arises from the behavior of the integrand: when $\aaa(x)=0$ the integrand reduces to $\snr{Dw}^p$ ($p$-phase) and therefore has $p$-growth, whereas for 
$
\aaa(x)>0$, it grows like $\snr{Dw}^q$ ($q$-phase). These properties make structures like $\mathcal D_{\aaa}$ suitable for modeling highly anisotropic materials and composites. Defining $\mathcal H(x,z):= |z|^p+\aaa(x)|z|^{q}$ for $z \in \er^n$, it follows that 
$
\mathcal{R}_{\mathcal{H}}(x, z) \leq c (p,q), 
$
so that  $\mathcal D_{\aaa}$ is still uniformly elliptic. Moreover, when $\mf{a}=\aaa_0\in \er$ is a constant function, more advanced regularity theory ensures that minimizers of 
$\mathcal D_{\aaa_0}$ are $C^{1, \beta}$-regular for some $\beta\equiv \beta(n,p,q)>0$. It is then natural to expect, in view of the available perturbation arguments for uniformly elliptic integrands, that when $\aaa$ is H\"older continuous, minimizers of 
$\mathcal D_{\aaa}$ should have H\"older continuous gradients. Surprisingly, this is not the case. Starting from \cite{sharp, vassili}, for every choice of $\alpha \in (0,1]$ and $\eps>0$, \cite{irene2} provides examples where 
$ \aaa \in C^{0,\alpha}$, $p<n<n+\alpha <q$ 
and there exist minimizers of $\mathcal D_{\aaa}$ whose set of essential discontinuities (the singular set) is a fractal with Hausdorff dimension exceeding $n-p-\eps$. Moreover, there exist minimizers that do not belong to $W^{1,q}_{\loc}$.  Furthermore, in \cite{anna} another fractal example is constructed in which the singular set has Hausdorff dimension $n-p$, the maximal dimension allowed for $W^{1,p}$-functions, and the minimizer fails to belong to any Sobolev space better than $W^{1,p}$. In such examples one can replace the integrand $\mathcal H(x,z)$  with its non-degenerate version  $(x,z)\mapsto H_1(z)^{p/2}+\mf{a}(x)H_1(z)^{q/2}$, which is smooth with respect to the gradient variable. All in all, we observe  minimizers of a uniformly elliptic, scalar functional, with a smooth integrand with respect to gradient variable and H\"older coefficients, that are just as bad as any other competitor. In the examples in  \cite{sharp, irene2} we also see that 
\eqn{pqno}
$$
\frac qp > 1+\frac \alpha n
$$
and that $q/p$ can be taken arbitrarily close to $1+\alpha/n$. This is a key point. In fact, after the first regularity results in \cite{sharp}, it was later proved by Paolo Baroni, Maria Colombo and the second author of this paper that minimizers of $\mathcal D_{\aaa}$ have H\"older continuous gradient when $\aaa \in C^{0,\alpha}(\Omega)$ and $q/p \leq  1 + \alpha/n.$ 

It is worth remarking that, while functionals like $\mathcal D_{\aaa}$ are uniformly elliptic for every choice of $\aaa$, they may still exhibit a weaker form of nonuniform ellipticity. To frame this, we introduce the concept of nonlocal ellipticity ratio for a general functional $\mathcal{F}$ in \rif{basicF}
$$
\mf{R}_{F}(z,\BBB) := \frac{\sup_{x\in \BBB}\Lambda(x,z)}{\inf_{x\in \BBB}\lambda(x,z)}
$$
for $\snr{z}\not =0$ and for every ball $\BBB\subset \Omega$. Note that we interpret $\mf{R}_{F}(z,\BBB)=\infty$ if the denominator vanishes. 
It follows that $\mathcal R_{F}(x,z)\leq \mf{R}_{F}(z,\BBB)$ whenever $(x,z) \in \BBB\times  (\er^n\setminus\{0_{\er^n}\})$. We shall say that the functional $\mathcal{F}$ is softly nonuniformly elliptic if it is uniformly elliptic but there exists a ball $\BBB\subset \Omega$ such that $z \mapsto \mf{R}_{F}(z,\BBB)$ is not bounded. This is the case for the functional $\mathcal D_{\aaa}$ whenever $\aaa $ vanishes but is not identically zero. Taking for instance $\aaa(x)=\snr{x}^\alpha$, we find 
$ \mf{R}_{\mathcal D_{\aaa}}(z,B_r(0))\approx 1 + r^{\alpha}\snr{z}^{q-p}$. The competition between growth with respect to $\snr{z}$ and decay of $r$ at zero in keeping $\mf{R}_{\mathcal D_{\aaa}}$ bounded of this last example is a key to understand the subtle interaction between gradient and coefficients in nonuniformly  elliptic problems encoded in condition \rif{pqno} (indeed see Theorem \ref{main} below).  

As demonstrated in \cite{sharp}, examples of irregular minima arise from the so-called Lavrentiev phenomenon, first identified by Zhikov for the functional $\mathcal{D}_{\aaa}$ \cite{vassili}. This occurs when minimizing a functional over a dense subspace yields a strictly larger value than the minimum over the entire space. In our situation this means
\eqn{Lav}
$$
\scalebox{0.92}{
$ \displaystyle 
\inf_{w\in u_{0}+W^{1,p}_{0}(\BBB)}\mathcal{D}_{\aaa}(w,\BBB)<\inf_{w\in u_{0}+W^{1,q}_{0}(\BBB)}\mathcal{D}_{\aaa}(w,\BBB)$}
$$
provided $q>p$, where $\BBB=B_1(0)$ and $u_0\in W^{1,\infty}(\BBB)$. As a consequence, one can construct minimizers that do not belong to $W^{1,q}_{\loc}(\BBB)$. Let us briefly describe a construction for \rif{Lav} in the two-dimensional case $n=2$. We follow  \cite{sharp, vassili}, thus considering the range $p<2<2+\alpha <q$.  
There are three main characters in this story: the malicious competitor $u_*\in W^{1,p}(\BBB)$, the honest boundary datum $u_0\in W^{1,\infty}(\BBB)$, and the cheating coefficient $\aaa\in C^{0, \alpha}(\BBB)$, who flirts with both. In polar coordinates $x_1=\rr \cos \theta$, $x_2=\rr \sin \theta$, and  with $\texttt{h}>0$ to be chosen, these are
$$
\scalebox{0.89}{
$ \displaystyle 
u_{*}(\rr,\theta):=\texttt{h}\begin{cases}
    \displaystyle
\ \sin(2\theta)\quad &\mbox{if} \ \ 0\le \theta\le \pi/4\vspace{.1mm}\\
\displaystyle
\ 1\quad &\mbox{if} \ \ \pi/4< \theta\le 3\pi/4\vspace{.1mm}\\
\displaystyle
\ \sin (2\theta-\pi)\quad &\mbox{if} \ \ 3\pi/4< \theta\le 5\pi/4\vspace{.1mm}\\
\displaystyle
\ -1\quad &\mbox{if} \ \ 5\pi/4< \theta\le 7\pi/4\vspace{.1mm}\\
\ \sin(2\theta)\quad &\mbox{if} \ \ 7\pi/4< \theta < 2\pi,
\end{cases}$}
$$ 
$u_0(x):= \snr{x}^2u_{*}(x)$ and finally $
a(\rr,\theta):=\rr^{\alpha}\max\left\{-\cos(2\theta),0\right\}=\snr{x}^{\alpha-2} (x_2^2-x_1^2)_{+}$. 
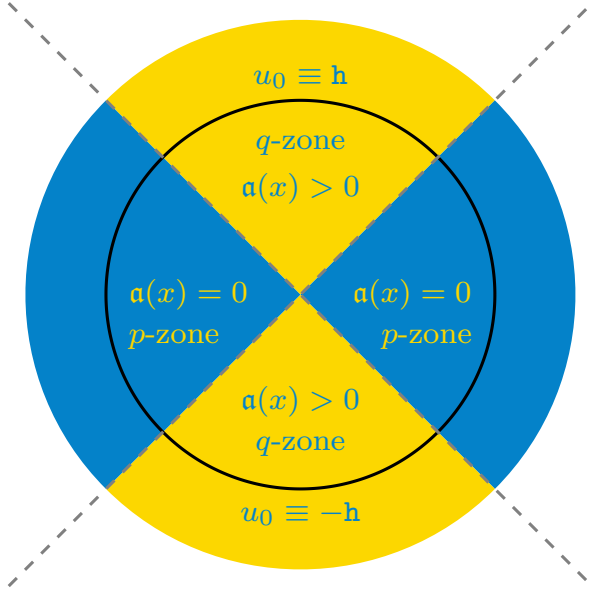
\begin{figure}[!ht]
    \centering
\begin{tikzpicture}[scale=1.45]
    \begin{axis}[
        hide axis,
        axis equal image,
        clip=false,
        width=8cm,
        height=8cm,
    ]
        \path[fill=gialloparma!100, draw=none] 
            (axis cs:0,0) -- (axis cs:1,1) arc (45:135:1.4142) -- cycle;
        
        \path[fill=bluparma!100, draw=none] 
            (axis cs:0,0) -- (axis cs:1,-1) arc (-45:45:1.4142) -- cycle;
        
        \path[fill=gialloparma!100, draw=none] 
            (axis cs:0,0) -- (axis cs:-1,-1) arc (-135:-45:1.4142) -- cycle;
        
        \path[fill=bluparma!100, draw=none] 
            (axis cs:0,0) -- (axis cs:-1,1) arc (135:225:1.4142) -- cycle;
        
        \addplot[
            domain=0:360,
            samples=500,
            thick,
            black 
        ] ({cos(x)}, {sin(x)});
        
        \addplot[
            domain=-1.5:1.5,
            thick,
            dashed,
            gray
        ] ({x}, {x});
        
        \addplot[
            domain=-1.5:1.5,
            thick,
            dashed,
            gray
        ] ({x}, {-x});

        \node[anchor=south, font=\footnotesize, bluparma!100!black] at (axis cs:0,0.65) {$q$-zone};
        \node[anchor=south, font=\footnotesize, bluparma!100!black] at (axis cs:0,0.4) {$\mathfrak{a}(x)>0$};
        \node[anchor=south, font=\footnotesize, bluparma!100!black] at (axis cs:0,-0.9) {$q$-zone};
        \node[anchor=north, font=\footnotesize, bluparma!100!black] at (axis cs:0,-0.4) {$\mathfrak{a}(x)>0$};
\node[anchor=south, font=\footnotesize, gialloparma!100!black] at (axis cs:0.65,-0.35) {$p$-zone};
        \node[anchor=west, font=\footnotesize, gialloparma!100!black] at (axis cs:0.2,0) {$\mathfrak{a}(x)=0$};
        \node[anchor=south, font=\footnotesize, gialloparma!100!black] at (axis cs:-0.65,-0.35) {$p$-zone};
        \node[anchor=east, font=\footnotesize,gialloparma!100!black] at (axis cs:-0.2,0) {$\mathfrak{a}(x)=0$};
        
        \node[anchor=north, font=\footnotesize, bluparma!100!black] at (axis cs:0,1.25) {$u_0\equiv\texttt{h}$};
        \node[anchor=north, font=\footnotesize, bluparma!100!black] at (axis cs:0,-1) {$u_0\equiv-\texttt{h}$};
        
    \end{axis}
\end{tikzpicture}
    \caption{Zhikov's checkerboard (in {\href{https://en.wikipedia.org/wiki/Parma}{\textcolor{bluparma}{Parma colours}}}). The yellow cone $\{x_2^2-x_1^2>0\}$ is the ``$q$-zone", where  $\aaa>0$ and you pay a lot. In the blue ``$p$-zone" $\{x_2^2-x_1^2\leq 0\}$, the coefficient $\aaa$ vanishes and you get a free lunch. }
\end{figure} Since $p<2$ and $\snr{Du_*}\lesssim 1/\snr{x}$, we have $u_{*}\in W^{1,p}(\BBB)$; similarly, $u_{0}\in W^{1,\infty}(\BBB)$. The key fact is that $a(x)\snr{Du_{*}}=0$ so that $\mathcal D_{\aaa}(u_{*};\BBB)=\nr{Du_*}_{L^p(\BBB)}^p= c_1(n,p) \texttt{h}^p$. By employing Direct Methods (see Section \ref{zerosec}) we find a unique minimizer $u\in u_{0}+W^{1,p}_{0}(\BBB)$ in the sense that the inf in the left-hand side of \rif{Lav} is attained at $w=u$. As $u_{*}\equiv u_{0}$ on $\partial \BBB$ it follows 
$$
\scalebox{0.92}{
$ \displaystyle 
\inf_{w\in u_{0}+W^{1,p}_{0}(\BBB)}\mathcal{D}_{\aaa}(w,\BBB)=\mathcal D_{\aaa}(u,\BBB)\le \mathcal D_{\aaa}(u_{*},\BBB) \leq c_1  \texttt{h}^p. $}
$$On the other hand, for any given $w\in u_{0}+C^{\infty}_{c}(\BBB)$ it is possible to show that there exists a constant $c_2\equiv c_2(n,q)\geq 1$ such that $\texttt{h}^{q}\leq c_2\mathcal D_{\aaa}(w;\BBB) $. By a standard density argument, this implies 
$$
\scalebox{0.92}{
$   \displaystyle 
\frac{\texttt{h}^q}{c_2}\leq \inf_{w\in u_{0}+W^{1,q}_{0}(\BBB)}\mathcal{D}_{\aaa}(w,\BBB) \,.$}
$$
Taking $\texttt{h}$ large enough to have $c_1c_2< \texttt{h}^{q-p}$ we reach \rif{Lav}. The key takeaway is that while a $W^{1,q}$-function, being H\"older continuous\footnote{Morrey's embedding claims that $W^{1,q}\subset C^{0, 1-n/q}$ when $q>n$. Here we are taking $q>2=n$.}, cannot avoid the yellow $q$-zone where movement is energetically costly, a  $W^{1,p}$-function for $p<2$ can jump and avoid that zone, concentrating the oscillations in the $p$-zone, thus minimizing energy expenditure. A singularity forms at the origin provided $\texttt{h}$ is sufficiently large, meaning that the boundary datum forces large oscillations in the $q$-zone for $W^{1,q}$-competitors. The example we have described features a cone-shaped coefficient $\aaa$ with vertex at the origin, generating in fact a one-point singularity of the minimizer. By distributing such cone shapes  over a Cantor set via more complex coefficients $\aaa$, more minimizers with large, fractal singular sets can be constructed, as shown in \cite{anna, irene2}.  

In the next sections we shall explore some classical and new ways of proving gradient bounds for nonlinear equations, eventually arriving at the proof of Schauder estimates in the nonuniformly elliptic setting in Section \ref{stimesec}. The Lavrentiev phenomenon will come back in Section \ref{finale}. 
 
\section{Inside De Giorgi's approach}\label{inside}
De Giorgi's celebrated approach to regularity of $W^{1,2}$-solutions to linear equations with bounded and measurable coefficients \rif{eq1} avoids any perturbation argument. As such, it does not rely on the linearity of the equations considered and, in fact, it can be extended to nonlinear equations. It is based on certain Caccioppoli type inequalities for truncations of solutions, namely 
\eqn{caccioppoli} 
$$
\scalebox{.95}{
$ \displaystyle 
   \int_{B_{\rr/2}} |D(v-\kk)_+|^2 \dx  \lesssim_{n,L}  \frac{1}{\rr^2} \int_{B_{\rr}}(v-\kk)_+^2\dx$}\,.$$
This holds for any choice of concentric balls $B_{\rr/2}\Subset B_{\rr}\Subset \Omega$ and numbers $\kk \in \er$. The proof of \rif{caccioppoli} is straightforward: one multiplies \rif{eq1}$_1$ by $\varphi= (v-\kk)_+\eta^2$, where $\eta\in C^{\infty}$, $\mathds{1}_{B_{\rr/2}}\leq \eta \leq \mathds{1}_{B_{\rr}}$, $\nr{D\eta}_{L^{\infty}(B_{\rr})}\lesssim 1/\rr$. Then \rif{caccioppoli} follows integrating by parts and using ellipticity \rif{eq1}$_2$ and Young’s inequality.
 The inequality \rif{caccioppoli} alone implies the local boundedness of $v$. The proof uses a nonlinear iteration based on a combination of \rif{caccioppoli} with Sobolev embedding theorem, and ultimately leading to  
  \eqn{reversa}
$$
 \nr{v_+}_{L^\infty(\BBB/2)}\lesssim_{n,L}\  \avnorm{v_+}_{L^2(\BBB)}
$$ 
whenever $\BBB \Subset \Omega$ is a ball. 
Applying \rif{reversa} to $-v$, that still solves \rif{eq1}, leads to replace $v_+$ by $v$ in \rif{reversa} thereby concluding  with 
the same $L^2-L^\infty$ local estimate valid for harmonic functions. 
This phenomenon is characteristic of such methods: despite their highly nonlinear nature they often recover the same estimates that hold for linear equations, such as the Laplace equation $\Delta u =0$. This is exactly the viewpoint of what is nowadays called Nonlinear Potential Theory: reproducing, in the nonlinear setting, the fine properties of solutions known in the linear one. 
Estimates of the type in \rif{caccioppoli} also imply local H\"older continuity of $u$, but here we are mainly interested in $L^{\infty}$-estimates. Inequality  \rif{caccioppoli}, and therefore \rif{reversa}, actually holds for so-called $W^{1,2}$-regular subsolutions, that is functions $v\in W^{1,2}_{\loc}(\Omega)$ satisfying 
\eqn{sottosola}
$$
\scalebox{0.9}{
$ \displaystyle -\diver(a(x)Dv)\leq 0 \Longleftrightarrow \int_{\Omega} a_{ij}(x)D_{j}vD_{i}\varphi\dx\leq 0 $}$$
whenever $0\leq \varphi\in C^{\infty}_0(\Omega)$. This paves the way for gradient estimates for uniformly elliptic nonlinear equations like \rif{el0}. The key is the so-called Bernstein method, which, roughly speaking, asserts that certain convex functions of $\snr{Du}$ are subsolutions to linear elliptic equations as in \rif{eq1}. To illustrate this in a model situation, we consider assumptions \rif{assf}$_{1-4}$, in the uniformly elliptic setting $p=q$ and when $F(x, z)\equiv F(z)$, and outline a proof of local Lipschitz continuity of  $W^{1,p}$-solutions to \rif{el0}. This covers the case of the $p$-Laplace equation $\diver (\snr{Du}^{p-2}Du)=0.$ The whole point is to show that 
 \eqn{lav}
$$
v := (\snr{Du}^{2}+\mu^{2})^{p/2}
$$
solves 
\rif{sottosola} with 
\eqn{linearized}
$$
 a_{ij}(x):= (\snr{Du(x)}^{2}+\mu^{2})^{\frac{2-p}{2}}\partial_{z_{i}z_{j}}F(Du(x)).
$$
We provide a formal proof in the non-degenerate case $\mu>0$, meaning that we demonstrate a priori estimates for more regular, say $C^1$-regular, solutions $u$ with $D_ju \in W^{1,2}$ for every $j$; all the computation can be then rigorously justified via suitable approximations methods that also allow to handle the case $\mu=0$. Differentiating \rif{el0} in the $s$ direction, for $s\in \{1, \ldots, n\}$, that is, choosing $\varphi\equiv D_s\varphi$ in \rif{debole} and integrating by parts, we obtain
\eqn{diffe}
$$
\int_{\Omega} \partial_{z_{i}z_{j}}F(Du)D_{js}u D_{i}\varphi\dx=0
$$
again for every $\varphi \in C^{\infty}_0(\Omega)$. 
 This step is already nontrivial, requiring careful use of difference quotient techniques due to the degeneracy of the equation (indeed, here we are taking $\mu>0$ but estimates are independent of $\mu$).
In \rif{diffe} we then choose $\varphi=  \eta D_su$ with $0\leq \eta\in C^{\infty}_0(\Omega)$ and sum over $s$.  
Using that $D_{j}v= p(\snr{Du}^{2}+\mu^{2})^{(p-2)/2} D_{sj}uD_su$,  we get (Einstein notation over  repeated indexes $i,j,s$ here)
$$
\scalebox{0.88}{
$ \displaystyle \int_{\Omega} a_{ij}(x)D_{j}v  D_{i}\eta \dx
+p\int_{\Omega} \eta  \partial_{z_{i}z_{j}}F(Du)D_{js}u D_{is}u \dx=0.$}
$$
The second integral is non-negative by 
\rif{assf}$_3$ and $\eta\geq 0$, which implies \rif{sottosola}. Moreover, by \rif{assf}$_{3,4}$ with $p=q$, the matrix $a(\cdot)$ in \rif{linearized} satisfies upper and lower bounds as in \rif{eq1} so that we can apply \rif{reversa} to $v$. Recalling \rif{lav} this finally translates \rif{reversa}  into
$ \nr{Du}_{L^\infty(\BBB/2)}\lesssim \ 
 \avnorm{Du}_{L^{2p}(\BBB)}+\mu$, that, after a standard interpolation/covering argument, gives 
the desired $L^{p}-L^\infty$ local gradient bound 
\eqn{reversa2}
$$
 \nr{Du}_{L^\infty(\BBB/2)}\lesssim \ 
 \avnorm{Du}_{L^{p}(\BBB)} + \mu
$$
with the implied constant depending only on $n,p,L$.

\section{Nonlinear potentials
}\label{dgsec} 
De Giorgi's iteration can be framed within the context of Nonlinear Potential Theory. This perspective  finds its roots in Kilpel\"ainen and Malý's seminal paper \cite{kilp} and has been extensively employed in gradient regularity theory, beginning by \cite{min}; see \cite{cris1, kumi} for an overview and references. In the foundational paper \cite{mazya}, Havin and Maz'ya introduced a class of nonlinear potentials designed for studying nonlinear,  possibly degenerate equations, in a manner analogous to the role of Riesz potentials in linear theory. Let $\sigma>0$, $\vartheta\geq 0$ and $f\in L^{1}(B_{r}(x))$; the nonlinear Havin-Maz'ya potential ${\bf P}_{\sigma}^{\theta}(f;\cdot)$ is defined as
$$
{\bf P}_{\sigma}^{\vartheta}(f;x,r) := \int_0^r \varrho^{\sigma} 
\avnorm{f}_{L^1(B_{\varrho}(x))}^{\vartheta}
\frac{\d\varrho}{\varrho}.
$$
The main feature of ${\bf P}_{\sigma}^{\vartheta}$ is in a sense its ability to capture the behavior of $f$ at each scale
\eqn{scale}
$$
\sum_{k=1}^\infty  (r/2^{k})^\sigma 
\avnorm{f}_{L^1(B_{r/2^k}(x))}^{\vartheta}\lesssim {\bf P}_{\sigma}^{\vartheta}(f;x,r) \,.
$$
For suitable choice of $\sigma, \vartheta$ the potentials ${\bf P}_{\sigma}^{\vartheta}$ recover classical Riesz or Wolff potentials. For instance, when $\sigma=\vartheta=1$ we find the usual (truncated) Riesz potential ${\bf I}_{1}^{f}$
\eqn{potenziali}
$$
\scalebox{0.97}{
 $ \displaystyle 
{\bf P}_{1}^{1}(f;x,r)\approx {\bf I}_{1}^{|f|}(x,r):= \int_0^r \frac{\nr{f}_{L^1(B_\rr(x))}}{\rr^{n-1}}\frac{\d\varrho}{\varrho}$}.
$$
Note that a simple annuli decomposition and \rif{scale}-\rif{potenziali} give
$$
\scalebox{1}{ $ \displaystyle 
\int_{B_{r/2}(x)} \frac{\snr{f(y)}}{\snr{y-x}^{n-1}}\dy \lesssim {\bf I}_{1}^{f}(x,r)\,.$}
$$ 
We know for which functions ${\bf P}_{\sigma}^{\vartheta}$ is bounded. 
\begin{thm}\label{crit} 
Let $B_{\tau}\Subset B_{\tau+r}\subset \mathbb{R}^{n}$, $n\geq 2$, be two concentric balls with $\tau, r\leq 1$, $f\in L^{1}(B_{\tau+r})$ and let $\sigma,\vartheta>0$ be such that $ n\vartheta>\sigma$. Then
 \eqn{stimazza0}
$$
 \scalebox{0.95}{
 $ \nr{{\bf P}_{\sigma}^{\vartheta}(f;\cdot,r)}_{L^{\infty}(B_{\tau})} \leq c\, \|f\|_{n\vartheta/\sigma, \vartheta;B_{\tau+r}}^{\vartheta} $}$$
holds with $c\equiv c (n,\vartheta,\sigma)$. In particular, 
 \eqn{stimazza}
$$
\scalebox{0.95}{
 $ \nr{{\bf P}_{\sigma}^{\vartheta}(f;\cdot,r)}_{L^{\infty}(B_{\tau})} \leq c\, \|f\|_{L^{m}(B_{\tau+r})}^{\vartheta} $}$$
holds for $m > n\vartheta/\sigma$ and with $c$ depending on $m$ too. 
\end{thm}
A proof can be found in \cite{cris3}. In \rif{stimazza0} the Lorentz space $L (n\vartheta/\sigma, \vartheta)$ appears. A measurable function $w$ belongs to $L(s,\gamma)(\Omega)$, $s,\gamma\in (0,\infty)$, when the quantity
$$
 \scalebox{0.9}{
 $ \displaystyle
\|w\|_{s, \gamma;\Omega} := \left(s\int_0^\infty (\lambda^s|\{x \in \Omega  :  |w(x)|> \lambda\}|)^{\gamma/s}\, \frac{d\lambda}{\lambda}\right)^{1/\gamma}$}
$$
is finite. 
When instead $s >0$ and $\gamma =\infty$ these are the usual Marcinkiewicz spaces\footnote{For $w$ to belong to $L(s,\infty)(\Omega)$ it is indeed required that 
$
\|w\|_{s, \infty;\Omega} := \sup_{\lambda >0} \lambda |\{x \in \Omega  :  |w(x)|> \lambda\}|^{1/s} < \infty.$}. 
In Lorentz spaces, the former index dominates in the sense that, when $\Omega\subset \mathbb{R}^{n}$ has finite measure, we have the following continuous inclusions:
$$ 
\begin{cases}
\,  \mbox{$L(s,s)(\Omega)=L^{s}(\Omega)$ for all $s >0$},\\
\, L(s_{1}, \gamma_{1})(\Omega)\subset L(s_{2}, \gamma_{2})(\Omega)\\
 \, \qquad  \mbox{ for all $0<s_{2}<s_{1}<\infty$, $\gamma_{1},\gamma_{2}\in (0,\infty]$,}\\
\, L(s,\gamma_{1})(\Omega)\subset L(s,\gamma_{2})(\Omega)\\
\, \qquad \mbox{ for all $s\in (0,\infty)$, $0<\gamma_{1}\leq \gamma_{2}\le\infty$}.
 \end{cases}
$$
The first two lines in the above display allow to derive \rif{stimazza} from \rif{stimazza0}. Lorentz spaces are particularly useful for asserting borderline regularity results for solutions to elliptic and parabolic equations, as we shall see at the end of this section. We now turn to the promised nonlinear potential-theoretic version of De Giorgi's iteration.
\begin{thm}\label{revlem}
Let $B_{r}(x)\subset \mathbb{R}^{n}$ be a ball, $n\ge 2$, $f \in L^1(B_{2r}(x))$, and constants $\chi >1$, $\sigma, \vartheta,\ti{c},M_{0}>0$ and $\kappa_0, M_{*}\geq 0$. Assume that $v \in L^2(B_{r}(x))$ is such that for all $\kk\ge \kk_{0}$  and $\rr \leq r$ the inequality
\begin{flalign}
& \displaystyle\avnorm{(v-\kk)_{+}}_{L^{2\chi}(B_{\rr/2}(x))}   \label{revva} \\ & \, \le \ti{c}M_{0}\avnorm{(v-\kk)_{+}}_{L^{2}(B_{\rr}(x))} \notag
+\ti{c} M_{*}\rr^{\sigma}\avnorm{f}_{L^{1}(B_{\rr}(x))}^{\vartheta}
\end{flalign}
holds. If $v$ has a precise representative at $x$ in the sense that 
$$
 \scalebox{0.97}{$\displaystyle v(x) = \lim_{\rr\to 0} \frac{1}{\snr{B_{\rr}(x)}}\int_{B_{\rr}(x)}v(y)\dy\,,$}$$ 
then
\begin{flalign}
v(x)  & \displaystyle \le \kk_{0}+cM_{0}^{\frac{\chi}{\chi-1}}\avnorm{(v-\kk_0)_{+}}_{L^{2}(B_{r}(x))}\notag\\
 & \quad \  +cM_{0}^{\frac{1}{\chi-1}} M_{*}\mathbf{P}^{\vartheta}_{\sigma}(f;x,2r)\label{stima}
\end{flalign}
holds with $c\equiv c(n,\chi,\sigma,\vartheta,\ti{c})$.  
\end{thm}
The proof, presented in \cite{cris3}, builds upon a refined version of De Giorgi’s iteration, employing a recursive selection of levels $\kappa$ - an idea originally introduced in \cite{kilp}.  Indeed, let us see how to get estimate \rif{reversa} from \rif{caccioppoli} via Theorem 3. By using Sobolev embedding theorem and rescaling, estimate \rif{caccioppoli} leads to the  reverse H\"older type inequality:
$$
 \scalebox{0.92}{
 $ \displaystyle \avnorm{(v-\kk)_{+}}_{L^{2\chi}(B_{\rr/2}(x))}
\lesssim_{n,L}\, 
   \avnorm{(v-\kk)_{+}}_{L^{2}(B_{\rr}(x))}
  $ }$$
for every ball $B_{\rr}(x) \Subset \Omega $ and  $\kk \in \er$, where $\chi := n/(n-2)>1$ for $n>2$, and any $\chi > 1$ works when $n=2$. Theorem \ref{revlem} applied with $\kappa_0=0$ implies that
\eqn{stimina}
$$
v(x)\lesssim_{n,L}\, 
\avnorm{v_+}_{L^{2}(B_{r}(x))}   
$$
holds whenever  $B_{r}(x)\Subset \Omega$ is a ball and $v$ 
 has the precise representative at $x$ - an occurrence that holds almost everywhere. Finally, consider a ball $\BBB \Subset \Omega$ with radius $\tau$. Applying \rif{stimina} with $r=\tau/4$ for (almost) every $x\in \BBB/2$ then easily leads to \rif{reversa}. At the core of De Giorgi's iteration lie self-improving effects of the type seen in \rif{reversa}: a suitable choice of the levels $\kappa$ allows to implement a geometric iteration where the improvement in the exponent, from $2$ to $2\chi$, accelerate convergence and leads to a pointwise bound for $v$. Two key aspects of Theorem \ref{revlem} are worth highlighting. First, the terms involving $f$ in \rif{revva}; these are later encoded in the potential $\mathbf{P}^{\vartheta}_{\sigma}(f;\cdot)$ in \rif{stima}. Second, the precise tracking of the constants $M_{0},M_{*}$ from \rif{revva} to \rif{stima}, which is usually overlooked in regularity proofs but proves essential when applying Theorem \ref{revlem} to nonuniformly elliptic problems, as we shall see in Section \ref{stimesec}. Regarding the first aspect, we recall a result from \cite{kumi} demonstrating the effectiveness of potentials in obtaining sharp gradient bounds in nonlinear, uniformly elliptic problems. Consider a vector valued solution $u$ to the $p$-Laplace system $\diver (\snr{Du}^{p-2}Du)=f$, $p\geq 2$. Surprisingly, $Du$ admits a pointwise estimate through linear potentials, mirroring the one available for the Poisson equation: 
\eqn{pppot}
$$
 \snr{Du(x)}^{p-1} \lesssim_{n,p} {\bf I}_{1}^{\snr{f}}(x,\rr)  
 +  \avnorm{Du}_{L^{1}(B_{\rr}(x))}^{p-1}\,.
 $$
See \cite{kumi} for the precise notion of solution for which \rif{pppot} holds. In particular, \eqref{pppot} remains valid when 
$f$ is a vector-valued Radon measure, in which case solutions are not generally $W^{1,p}$-regular. 
Estimate \rif{pppot} actually holds for general quasilinear equations (but not for general systems!) of $p$-Laplace type and allows to reduce local gradient bounds of solutions to nonlinear equations to those for the Poisson equation. For instance, it immediately leads to a sharp locally Lipschitz continuity criterion for solutions to the $p$-Laplace system, which coincides with the one known for $\Delta u=f$, namely $f\in L(n,1)$. Remarkably, such a condition is independent of $p$, confirming an old conjecture usually attributed to Ural'tseva. For this see also former work of Duzaar and the second author of this paper. In the linear case, this follows from a famous result of Stein, which states that $Dw\in L(n,1)$ implies the continuity of $w$. In fact, in \cite{kumi} it is proved that $f\in L(n,1)$ ensures the continuity of $Du$ for the $p$-Laplace system.

 \section{Schauder estimates}\label{stimesec}
  The problem of deriving Schauder estimates in the nonuniformly elliptic case has been successfully tackled in \cite{cris3, cris4}, where some new approaches to gradient bounds have been developed. For an initial exposition, we shall restrict ourselves to the simpler case of functionals of the type 
  \eqn{modello}
$$
 \mathcal G(w, \Omega ):=\int_{\Omega}\ccc(x)G(Dw) \dx\,, \ \ 1 \leq \ccc \leq L
$$
for which results are more easily displayed. Nevertheless, from the viewpoint of purely a priori regularity estimates, this model already encapsulates all the difficulties arising in the general case \rif{basicF}. The assumptions on $G$ we consider are those in \rif{assf}, adapted to the case where there is no explicit $x$-dependence in the integrand. The coefficient $\ccc$ instead satisfies
 \eqn{coeff}
$$
|\ccc(x)-\ccc(y)| \leq L|x-y|^{\alpha},  \  \  \alpha \in (0,1]
$$
whenever $x, y\in \Omega$. 
\begin{thm}[Nonuniformly elliptic Schauder]\label{main}
Let $u\in W^{1,1} (\Omega)$ be a minimizer of the functional $ \mathcal{G}$, under assumptions \eqref{assf}$_{1-4}$ (with $F\equiv G$), \eqref{coeff} and 
\eqn{pq}
$$
\frac qp < 1+\frac \alpha n.
$$ Then, $Du$ is locally H\"older continuous in $\Omega$. If $\mu>0$, then $u\in C^{1,\alpha}_{\loc}(\Omega)$ when $\alpha <1$. \end{thm}
The non-degeneracy assumption $\mu >0$ is necessary to obtain the sharp result $Du \in C^{0, \alpha}_{\loc}(\Omega;\er^n)$. Otherwise, well-known counterexamples - already valid in the case of the homogeneous $p$-Laplace equation - demonstrate that achieving H\"older continuity for every $\alpha <1$ is impossible. The main assumption in Theorem \ref{main} is \rif{pq}, which, apart from the equality case, complements the condition under which counterexamples to regularity arise in the nonautonomous case \rif{pqno}. See also Theorem \ref{mainr} and Corollary \ref{mainrr} in the next section.

The core of Theorem \ref{main} lies in establishing the local boundedness of $Du$. Once this is achieved, one can, with careful handling of some delicate details, adapt perturbation arguments typically used in the uniformly elliptic setting. Traditionally, gradient $L^\infty$-bounds under H\"older-continuous coefficients follow as a consequence of $C^{1,\alpha}$-bounds, which are, in turn, obtained via perturbation methods, i.e., Schauder estimates. However, due to the severe inhomogeneity of the available a priori estimates, this approach is not viable in the nonuniformly elliptic case, as already observed in Section \ref{gradb}. In \cite{cris3,cris4}, a method has been devised that, to the best of our knowledge, allows for the first time a direct derivation of $L^\infty$-gradient bounds without any quantitative appeal to $C^{1,\alpha}$-estimates. In a traditional approach, as outlined in Section \ref{inside}, one would differentiate the Euler-Lagrange equation 
\eqn{el}
$$
-\diver \left(\ccc(x) \partial_zG(Du)\right)=0
$$
to eventually prove that a function $v$ like in \rif{lav} is a subsolution to a linear elliptic equation as in \rif{eq1}. This is obviously not possible in this setting due to the presence of a non-differentiable coefficient $\ccc$ (that's why perturbation theory is used in such cases!). Instead, the idea is to interpret the H\"older continuity of $\ccc$ in terms of fractional differentiability and to differentiate \rif{el} in a fractional sense. We then discover that $v$, rather than satisfying a standard differential equation, verifies Caccioppoli-type estimates as in \rif{caccioppoli}, with the left-hand side replaced by a suitable fractional norm. Once the appropriate fractional Caccioppoli inequality is established, we can apply the nonlinear potential theoretic machinery described in Section \ref{dgsec}. This approach leads to $L^\infty$-gradient bounds \cite{cris3} and, when combined with additional higher integrability results, leads to the same bounds but under the optimal assumption \rif{pq} \cite{cris4}. We are now going to sketch the backbone of the proof of Theorem \ref{main}, warning the reader that, as usual in such cases, we shall only deal with a priori estimates for more regular minimizers, say $C^1$. This is eventually fixed via approximation. We shall confine ourselves to $L^\infty$-gradient bounds. Specifically,
\eqn{m.1}
$$
\scalebox{0.95}{
$ \displaystyle
 \nr{Du}_{L^\infty(\BBB/4)}\leq c \, 
 \avnorm{G(Du)}_{L^{1}(\BBB)} ^{\mf{b}}+c
$}
$$
holds with $c\equiv c(n,p,q,\alpha,L)\geq1$, $\mf{b}\equiv \mf{b}(n,p,q,\alpha)\geq 1$, for every ball  $\BBB\Subset \Omega$ with radius not exceeding $1$. To proceed, we recall that a map $w\colon \Omega \mapsto \er^n$ belongs to $W^{s,m}(\Omega;\er^n)$, $0<s<1$, $m\geq 1$, if $w\in L^m(\Omega;\er^n)$ and
$$
\scalebox{0.95}{
 $  \displaystyle [w]_{\tta,\mmm;\Omega}:= \left(\int_{\Omega} \int_{\Omega}  
\frac{|w(x)
- w(y) |^{\mmm}}{|x-y|^{n+\tta \mmm}} \dx \dy \right)^{1/\mmm}< \infty$}.
$$
Moreover, let $\texttt{t}>0$, $h \in \mathbb{R}^n$, and set $\Omega_{\texttt{t}}:=\left\{x\in \Omega\colon \dist(x,\partial \Omega)>\texttt{t}\right\}$;  the finite difference operators $\tau_{h}\colon L^{1}(\Omega;\mathbb{R}^{n})\mapsto L^{1}(\Omega_{|h|};\mathbb{R}^{n})$, $\tau^{2}_{h}\colon L^{1}(\Omega;\mathbb{R}^{n})\mapsto L^{1}(\Omega_{2|h|};\mathbb{R}^{n})$ are defined as $ \tau_{h}w(x):=w(x+h)-w(x)$ and  $\tau_{h}^{2}w:=\tau_{h}(\tau_{h}w)$, respectively. The proof of \rif{m.1} proceeds in three steps. 
 
{\em Step 1: Higher integrability of $Du$}. This means \eqn{almos}
$$
Du \in L^{\mf{q}}_{\loc}(\Omega;\er^n) \quad \mbox{for every $\mf{q}<\infty$}.
$$
The proof relies on Besov spaces techniques along an iteration scheme that, in a sense, resembles  Moser's iteration but with a linear rather than geometric gain in the exponents at each step. 
Indeed, one obtains estimates of the type 
\eqn{derive}
$$
\nr{\tau_{h}^{2}u}_{L^{\mf{q}_{i}}} \leq c_{i}\snr{h}^{t_{i}}
$$
valid for $h \in \er^n$ and $|h|$ suitably small. In \rif{derive} there appear three sequences $\{\mf{q}_{i}\}$,$\{t_{i}\}$,$\{c_{i}\}$. They are such that $
\mf{q}_{i}, c_{i} \to \infty$, $1< t_{i}  \to 1$. The key point here is that \rif{derive} implies that $Du \in L^{\mf{q}_{i}}_{\loc}(\Omega;\er^n)$. This follows by general embedding properties of Besov spaces via double difference operators, i.e., 
\eqn{thermic}
$$
 \sup_{0<\snr{h}\leq h_0}\snr{h}^{-t}\nr{\tau_{h}^{2}w}_{L^{\mf{q}}} <\infty\Longrightarrow \nr{Dw}_{L^{\mf{q}}}< \infty
$$
locally, where $h_0>0$, $t\in (1,2)$, $\mf q>1$. The result in \eqref{thermic} is closely related to the so-called thermic characterization of Besov spaces, as described in classical works by Stein and Triebel.  
 In \rif{derive} it is actually
$
\mf{q}_{i+1} = \mf{q}_{i} + \texttt{s}
$
for some $\texttt{s}\equiv \texttt{s}(n,p,q,\alpha)  >0$, thereby yielding a linear gain of integrability along the iteration (in classical Moser's iteration the gain is geometric). This is not fast enough to reach Lipschitz estimates already at this stage. 

{\em Step 2: A (renormalized) fractional Caccioppoli inequality}. 
 The function $v$ in \rif{lav} satisfies a fractional analogue of  \rif{caccioppoli}, i.e.,  
\begin{flalign}
&\scalebox{0.97}{$ \displaystyle  \nonumber  [(v-\kk)_{+}]_{\beta,2;B_{\rr/2}}^2 \le \frac{cM^{\mf{s}(q-p)}}{\rr^{2\beta}}\int_{B_{\rr}}(v-\kk)_{+}^{2}\dx$}\\ & \scalebox{0.97}{$ \displaystyle \qquad \qquad \ \ \,  +\frac{cM^{\mf{s}(q-p)}\rr^{2\alpha\kkk}}{\rr^{2\beta}}\int_{B_{\rr}}(\snr{Du}+1)^{\mf{p}}\dx$}\label{280}
\end{flalign}
holds for any $M\geq \max\{ \nr{Du}_{L^\infty(B_{\rr})}, 1\}$, concentric balls $B_{\rr/2}\Subset B_{\rr}\Subset \Omega$ with $\rr\leq 1$, and $\kappa \in \er$. In \rif{280}  the parameters $\beta, \kkk, \mf{p}$ satisfy $\beta \in (0,\alpha/(1+\alpha))$, $\kkk\in (0,1)$,  $\mf{p}> 1$  with a certain degree of flexibility in their choice, which involves a delicate interplay. Specifically, they will be chosen according to the following scheme:
\eqn{rule}
$$
\mbox{ $\beta\approx \frac{\alpha}{1+\alpha}\stackrel{\scalebox{.7}{requires}}{\Longrightarrow} \kkk \approx 1
\stackrel{\scalebox{.7}{requires}}{\Longrightarrow}\mf{p}$ is large}.
$$
Inequality \rif{280} ``improves" when taking $\beta$ and $\kkk$ close to their limit values in \rif{rule}, which is made possible by ``paying" a large 
$\mf{p}$, a trade-off eventually allowed by \rif{almos}.
Another technical point, and another basic difference with respect to \rif{caccioppoli}, is the presence of the multiplicative factor $M^{\mf{s}(q-p)}$ on the right-hand side of \rif{280}. The rationale behind this factor is tied to the use of iteration tools, such as those in Theorem \ref{revlem}, which require homogeneous estimates with respect to the iterating function $v$. However, such homogeneous estimates are not available for nonuniformly elliptic equations, presenting the same obstruction encountered with perturbation methods. To address this, the energy estimates must be formally renormalized, incurring the cost of  $M^{\mf{s}(q-p)}$, that eventually one hopes to reabsorb in some way. The outcome is \rif{280}, which is formally homogeneous in $v$ (but not in $\snr{Du}$ as a whole, as $M$ depends on $\nr{Du}_{L^\infty(B_\rr)}$). In this respect, note that in \rif{280} the number $\mf{s}$ is exactly the one in \rif{nuovino}. This highly parametrized form of estimates is peculiar in nonuniformly elliptic problems. Optimizing the arrangement of the parameters is crucial to obtaining sharp bounds such as \rif{pq}. The proof of \rif{280} is not straightforward. Typically, basic energy estimates as \rif{caccioppoli} involve the same order of differentiation of the operator and follow by plain testing. Here, the approach is nonlocal but applies to a local problem, which must, in a sense, be delocalized. This is accomplished through a domain decomposition technique that, from a conceptual viewpoint, resembles the atomic characterization of Besov functions in Littlewood–Paley theory. In a sense, the atoms used here are themselves solutions to autonomous nonlinear variational problems, such as those in Theorem \ref{mainaut}. This explains how the number 
$\sss$ in \rif{nuovino} emerges in \rif{280}. Arguments of this type were originally introduced by Kristensen and the second author of this paper to study singular sets in the Calculus of Variations via fractional Sobolev spaces.

{\em Step 3: Lipschitz regularity}. Using fractional Sobolev embedding in combination with \rif{280} yields
\begin{flalign}
& \notag \avnorm{(v-\kk)_{+}}_{L^{2\chi}(B_{\rr/2}(x))} \\
&\quad \label{28}  \le cM^{\mf{s}(q-p)/2}\avnorm{(v-\kk)_{+}}_{L^{2}(B_{\rr}(x))} \\ & \quad   \quad   +c M^{ \mf{s}(q-p)/2}\rr^{\alpha\kkk}\avnorm{(\snr{Du}+1)^{\mf p}}_{L^{1}(B_{\rr}(x))}^{1/2}\notag
\end{flalign}
whenever $B_{\rr}(x)\Subset \Omega$, $\rr\leq 1$, $\kk \in \er$, where this time $\chi = n/(n-2\beta)>1$. This enables the use of  Theorem \ref{revlem}. For the proof of \rif{m.1}, denote $\BBB\equiv B_{\tau}$, and take concentric balls $B_{\tau/4}\Subset  B_{\tau_{1}}\Subset   B_{\tau_{2}}\Subset B_{\tau/2}$ and set $r:=(\tau_{2}-\tau_{1})/16$ so that, if $x \in B_{\tau_1}$, then $B_{2r}(x)\Subset B_{\tau_{2}}$. We use \eqref{28} with $\rr\leq r$, $M:= \nr{Du}_{L^{\infty}(B_{\tau_{2}})}\geq \nr{Du}_{L^{\infty}(B_{\rr}(x))}$. Note that we can assume $M\geq1$ otherwise \rif{m.1} is trivial. Applying Theorem \ref{revlem} on $B_{r}(x)$ with $\kk_{0}=0$, $M_{0}= M_{*}= M^{\mf{s}(q-p)/2}$, $f= (\snr{Du}+1)^{\mf{p}}$, $\sigma= \alpha\kkk$, and $\vartheta= 1/2$, and yet recalling that $x\in B_{\tau_1}$ is arbitrary, we obtain
\begin{flalign*}
&\scalebox{0.88}{$ \displaystyle\nr{Du}_{L^{\infty}(B_{\tau_{1}})} \le\frac{c\nr{Du}_{L^{\infty}(B_{\tau_{2}})}^{\mf{t}}}{(\tau_{2}-\tau_{1})^{n/(2p)}}\nr{\snr{Du}+1}_{L^{2p}(\BBB/2)}$}\nonumber \\
&\scalebox{0.88}{$ \displaystyle\quad +c\nr{Du}_{L^{\infty}(B_{\tau_{2}})}^{\mf{t}}\nr{\mathbf{P}^{1/2}_{\alpha\kkk}\left((\snr{Du}+1)^\mf{p};\cdot,2r\right)}_{L^{\infty}(B_{\tau_{1}})}^{1/p}+c$}
\end{flalign*}
with $\mf{t}:= n(\mf{s}/4)(q/p-1)/\beta$. 
A key point is here: the bound on $q/p$ in \rif{pq} and the precise identity of the exponent $\sss$ in  \rif{nuovino}  allow to choose $\mf{p}, \kkk, \beta$ in \rif{rule} in such a way that $\mf{t}<1$. To bound the last term we appeal to \rif{stimazza} with $\sigma = \alpha \kkk$, $\vartheta =1/2$ and we take $m> n/(2\alpha \kkk)>1$ to get 
$$
\scalebox{0.88}{$
\nr{\mathbf{P}^{1/2}_{\alpha\kkk}\left((\snr{Du}+1)^\mf{p};\cdot,2r\right)}_{L^{\infty}(B_{\tau_{1}})}\lesssim \nr{\snr{Du}+1}_{L^{\mf{q}}(\BBB/2)}^{\mf{p}/2}$}
$$
for any $\mf{q}> n\mf{p}/(2\alpha \kkk)$. Thanks to \rif{almos} the right-hand side is under control and the last two displays yield 
\begin{flalign*}
&\nr{Du}_{L^{\infty}(B_{\tau_{1}})}\le\frac{c\nr{Du}_{L^{\infty}(B_{\tau_{2}})}^{\mf{t}}}{(\tau_{2}-\tau_{1})^{n/(2p)}}\nr{\snr{Du}+1}_{L^{2p}(\BBB/2)}\nonumber \\
&\qquad +c\nr{Du}_{L^{\infty}(B_{\tau_{2}})}^{\mf{t}}\nr{\snr{Du}+1}_{L^{\mf{q}}(\BBB/2)}^{\mf{p}/(2p)}+c.
\end{flalign*} 
Since $\mf{t}<1$ we use Young's inequality concluding with
\begin{flalign*}
& \notag \nr{Du}_{L^{\infty}(B_{\tau_{1}})} \le\scalebox{0.93}{$\frac{1}{2}$}\nr{Du}_{L^{\infty}(B_{\tau_{2}})}
\\ & \qquad +\frac{c\nr{\snr{Du}+1}_{L^{2p}(\BBB/2)}^{\frac{1}{1-\mf{t}}}}{(\tau_{2}-\tau_{1})^{\frac{n}{2p(1-\mf{t})}}}  +c\nr{Du}_{L^{\mf{q}}(\BBB/2)}^{\frac{\mf{p}}{2p(1-\mf{t})}}+c.
\end{flalign*}
The first term in the right-hand side cannot be directly reabsorbed into the left-hand side as $\tau_2>\tau_1$, but we can fix this via a well-known iteration lemma due to Giaquinta and Giusti\footnote{Let $\hhh\colon [r_1,r_2]\to \mathbb{R}$ be a non-negative and bounded function, and let $a,b, \gamma$ be non-negative numbers. Assume that the inequality 
$
\hhh(\tau_1)\le  \hhh(\tau_2)/2+a(\tau_2-\tau_1)^{-\gamma}+b
$
holds whenever $r_1\le \tau_1<\tau_2\le r_2$. Then 
$
\hhh(r_1)\lesssim_{\gamma} a(r_2-r_1)^{-\gamma}+b.
$}.
This allows us to ``reabsorb" that term anyway, thus obtaining 
$$
\scalebox{0.93}{$
\nr{Du}_{L^{\infty}(\BBB/4)}\le c\avnorm{Du}_{L^{2p}(\BBB/2)}^{\frac{1}{1-\mf{t}}}+c\nr{Du}_{L^{\mf{q}}(\BBB/2)}^{\frac{\mf{p}}{2p(1-\mf{t})}}+c$}.
$$
This and the bounds implicit in \rif{almos} allows us to conclude with \rif{m.1}; see \cite{cris4}. 

\section{A selection principle}\label{finale}
When transitioning from \rif{modello} to the general case \rif{basicF} the Lavrentiev phenomenon may arise, preventing regularity of minimizers. The strategy is then to consider a relaxed functional that discards singular minimizers caused by the Lavrentiev phenomenon and selects regular ones. From a variational perspective, such minimizers play a role analogous to that of energy solutions relative to general distributional solutions in the context of elliptic equations. Indeed, while energy solutions to \rif{eq1} - i.e., distributional solutions that are $W^{1,2}$-regular - are continuous, non-$W^{1,2}$ solutions can be unbounded, as originally shown by Serrin. Under assumptions \rif{assf}$_{1-4}$, the Lebesgue-Serrin-Marcellini relaxation of the functional $\mathcal{F}$ in \rif{basicF} is defined as
\begin{flalign}
&    \scalebox{0.88}{$\displaystyle 
\bar{\mathcal{F}}(w,U)$}  :=\label{relax}\\
 & \ \  \scalebox{0.88}{$ \displaystyle  \inf_{\{w_{i}\}\subset W^{1,q} (U)} \left\{ \liminf_{i} \mathcal{F}(w_{i},U)  \, \colon \,  w_{i} \deb w\  \mbox{in} \ W^{1,p}(U)  \right\} $}\notag
\end{flalign} 
whenever $w \in W^{1,p}(U)$ and $U\subset \Omega$ is an open subset. As $z \mapsto F(\cdot, z)$ is convex, De Giorgi-Ioffe Theorem implies that $\mathcal  F$ is weakly lower semicontinuous in $W^{1,p}$, ensuring $\mathcal{F}\leq \bar{\mathcal{F}}$. On the other hand, if $w \in W^{1,q}(U)$, then $\bar{\mathcal{F}}(w,U)=\mathcal{F}(w,U)$. Consequently, $\bar{\mathcal{F}}$ can be understood as a “non-pointwise” extension of $\mathcal{F}$ to $W^{1,p}$, assuming that $\mathcal{F}$ was originally defined only on $W^{1,q}$, that is the space where it is always finite due to \rif{assf}$_1$. With $\BBB\Subset \Omega$ being a ball, from \rif{relax} it follows that if  $\bar{\mathcal{F}}(w,{\BBB})<\infty $, then there exists a sequence $\{w_{i}\}\subset W^{1,q} (\BBB)$ such that
\eqn{appro}
$$\scalebox{0.95}{$\displaystyle \mbox{$w_{i} \deb w$ in $W^{1,p}(\BBB)$ and $ \mathcal{F}(w_{i},{\BBB}) \to \bar{\mathcal{F}}(w,{\BBB})$}$}.
$$ 
Finally, $w \mapsto \bar{\mathcal{F}}(w,\BBB)$ remains strictly convex, coercive and lower semicontinuous in the weak topology of $W^{1,p}(\BBB)$. A function $u \in W^{1,p}_{\loc}(\Omega)$ is a minimizer of $\bar{\mathcal{F}}$ if $\bar{\mathcal{F}}(u,{\BBB})<\infty $ and $\bar{\mathcal{F}}(u,{\BBB})\leq \bar{\mathcal{F}}(w,{\BBB})$ holds for every ball ${\BBB}\Subset \Omega$ and every $w \in u + W^{1,p}_0({\BBB})$. 
\begin{thm}[Selection of regular minimizers]\label{mainr}
Let $u\in W^{1,p}_{\loc}(\Omega)$ be a minimizer of $\bar{\mathcal{F}}$ under assumptions \eqref{assf} and \eqref{pq}. Then $u$ is also a minimizer of $\mathcal{F}$ and the regularity properties in Theorem \ref{main} hold. 
\end{thm}
The regularity part of Theorem \ref{mainr} goes as in Section \ref{stimesec} modulo an approximation argument based on \rif{appro}. Eventually, it is easy to see that local Lipschitz continuity implies that $u$ also minimizes $\mathcal{F}$ \cite{cris4}. 
\begin{cor}[Regularity via relaxation]\label{mainrr}
Let $u\in W^{1,p}_{\loc}(\Omega)$ be a minimizer of $\mathcal{F}$ under assumptions \eqref{assf} and \eqref{pq}, and also assume that 
\eqn{relax2}
$$\mathcal{F}(u, \BBB)=\bar{\mathcal{F}}(u, \BBB) \ \ \mbox{for every ball $\BBB\Subset \Omega$}.$$ Then $u$ is also a minimizer of $\bar{\mathcal{F}}$ and therefore enjoys the regularity properties of Theorem  \ref{main}. 
\end{cor} 
Indeed, using the minimality of $u$ and that $\mathcal{F}\leq \bar{\mathcal{F}}$, we have $\bar{\mathcal{F}}(u, \BBB) =\mathcal{F}(u, \BBB)\leq \mathcal{F}(w, \BBB)\leq  \bar{\mathcal{F}}(w, \BBB)$ whenever $w \in u+W^{1,p}_0(\BBB)$, so that $u$ also minimizes $\bar{\mathcal{F}}$ and Theorem \ref{mainr} applies. The functional $\bar{\mathcal{F}}$ was first introduced by Marcellini \cite{paolo} in the vectorial setting to provide a fascinating interpretation of cavitation phenomena in Nonlinear Elasticity through the appearance of the Lavrentiev phenomenon. Remarkable work of Fonseca and Malý \cite{irene} provides (weak) representation for $U\mapsto \bar{\mathcal{F}}(u, U)$ via a Radon measure $\mu_u$ such that $\mu_u(U)\leq  \bar{\mathcal{F}}(u, U)\leq \mu_u(\bar U)$ holds for every open subset $U\subset \Omega$, provided $ \bar{\mathcal{F}}(u, \Omega)$ is finite. This holds for general integrands $F(x,z)$ with $(p,q)$-growth conditions \rif{assf}$_2$ under bounds of the type \rif{central}, but not necessarily convex with respect to the gradient variable $z$. The singular part (with respect to the Lebesgue measure) of $\mu_u$ appears in connection to the Lavrentiev phenomenon. In this respect, the functional $w \mapsto \bar{\mathcal{F}}(w, \BBB)-\mathcal{F}(w, \BBB)\geq 0$ ($:=0$ when $\mathcal{F}(w, \BBB)=\infty$), called Lavrentiev gap, can be used to encode the presence of the Lavrentiev phenomenon in a ball $\BBB$. Indeed, in the setting of Corollary \ref{mainrr}, thanks to \rif{appro} it is possible to approximate $\mathcal{F}(u, \BBB)$ with a sequence $\{\mathcal{F}(u_i, \BBB)\}$, $\{u_{i}\}\subset W^{1,q} (\BBB)$, thereby excluding the Lavrentiev phenomenon. Summarizing, minimizers of $\bar{\mathcal{F}}$ are regular minimizers of $\mathcal{F}$. In the other direction, the regularity of minimizers of the ``pointwise" functional $\mathcal{F}$ is guaranteed provided the Lavrentiev phenomenon does not appear, which is obviously an essential condition for regularity. The question then arises whether the Lavrentiev phenomenon appears or not, or, more specifically in this setting, under which conditions  \rif{relax2} holds. This is the case for large classes of functionals, including for instance those as in \rif{modello}. See \cite{sharp}, where the point of view of Corollary \ref{mainrr} was introduced. An interesting twist is that in several cases condition \rif{pq} ensures \rif{relax2} independently of regularity estimates (and therefore can be combined with a priori estimates to prove regularity, as in Corollary \ref{mainrr}). This happens for instance with the double phase functional in \rif{doppio} and much larger classes of related integrals \cite{sharp}. We finally notice that a different, alternative definition to $\bar{\mathcal{F}}$, usually denoted by $\bar{\mathcal{F}}_{\loc}$ \cite{irene}, occurs when taking sequences  $\{w_{i}\}\subset W^{1,q}_{\loc} (U)$ in \rif{relax}; the results presented in this section continue to hold in that case. For this last functional, a stronger representation result holds \cite{irene}, namely, whenever $\bar{\mathcal{F}}_{\loc}(u,\Omega)$ is finite, there exists a Radon measure $\mu_u$ defined on $\Omega$ such that $\mu_u(U)=\bar{\mathcal{F}}_{\loc}(u,U)$ for all open subsets $U \subset \Omega$.

\vspace{2mm}
{\bf Acknowledgments.} 
This work was supported by the University of Parma through the action ``Bando di Ateneo 2024 per la ricerca".


       \end{document}